\input amstex
\documentstyle{amsppt}
\baselineskip 20pt
\magnification=1200
\NoBlackBoxes

\def\br#1{\langle #1 \rangle}
\def\ket#1{|#1 \rangle}
\def\Phit{\widetilde{\Phi}}

\centerline{ \bf Spinor representations of 
$U_q(\hat \frak {o}(N))$}
\vskip .2in 
\centerline{Jintai Ding}
\centerline{RIMS, Kyoto University}
\vskip .4in
\centerline{ abstract}
This is an extension of quantum spinor construction of 
$U_q(\hat {\frak gl}(n))$. We 
define quantum affine Clifford algebras based on the tensor category and 
the solutions of q-KZ equations, and  construct quantum spinor representations
of $U_q(\hat{\frak o}(N))$. 

\subheading{I. Introduction}

For affine Lie algebras, we have  quadratic constructions  of 
representations with Clifford, Weyl or Heisenberg algebras, which 
are  spinor or oscillator representations\cite{FF}\cite{KP}. 
One question is that if we can extend those constructions 
to the cases for quantum groups, a q-deformation 
of universal enveloping algebras of  Kac-Moody algebras
 discovered by Drinfeld \cite{D1} and Jimbo \cite{J1}. 
Lusztig \cite{L} showed the existence of 
such a q-deformation of the 
category of highest weight representations of Kac-Moody algebras for
generic  parameter q.
There also appeared many constructions of representations 
\cite{B}, \cite{FJ},  \cite{H}, \cite{J2}, etc.

In \cite{DF2} \cite{Di}, we proposed an invariant approach
to deform explicitly those classical quadratic 
constructions with  a q-analogue of matrix 
realization of classical Lie algebras. We managed to use such 
an  approach to 
define quantum Clifford and   Weyl algebras and quantum affine Clifford 
algebras   
using general representation theory of quantum groups and 
solutions of q-kz equations\cite{DO}. 
We  showed  that the explicit formulas for quantum Clifford and Weyl algebras
 match the ones actively studied in physics literature.  
Using 
those quantum algebras, we constructed 
 spinor and oscillator representations of 
quantum groups of classical types and 
spinor representations of $U_q(\hat {\frak{gl}(n)})$.
Representation theory allows  us to justify
that the quadratic
elements in quantum Clifford and Weyl algebras and quantum 
affine Clifford algebras provide the 
desired representations. In our case, an explicit verification
of Serre's relations for quantum groups is not necessary.
 The key idea   consists of reformulating 
 familiar classical constructions entirely
in terms of the  tensor category of highest weight representations and 
using Lusztig's result on q-deformation of this category and 
 solutions of q-kz equations to  define the
corresponding quantum structures. In the quantum case, the 
quasitriangular structure of the  tensor category introduced by Drinfeld
\cite{D3} plays the fundamental role.  We  would like to stress 
the central role of the universal Casimir operator and its inverse 
implied by 
the quantum structure\cite{Di}. 

In this paper, we  extend
 the idea in \cite{Di}  to the cases of  the 
spinor representations of 
quantum affine groups $U_q(\hat {\frak o}(N))$, which 
is simpler than that in \cite{Di} 
due to the structure of the corresponding modules. 
From what we have obtained, it is convincing that we can 
deform all the constructions given in \cite{FF} as 
straightforward extensions   of our constructions.

\cite{D1} \cite{J1}
Quantum group for $U_q(\frak {g})$ is defined as  an associative algebra 
generated by $e_i$, $f_i$ and $t_i$, $i=0,1..,n$, where $i$ corresponds 
to the index of the nodes of the Cartan matrix
(see Definition 1.1 in Section 1).  This algebra has 
a noncocommutative Hopf algebra structure 
 with comultiplication $\Delta$,
antipode $S$ and counit $\varepsilon$. 
 \cite{G} The affine
Kac-Moody algebra $\hat {\frak{g}}$ associated to a simple Lie algebra
$\frak{g}$ admits a natural realization as a central extension of the
corresponding loop algebra ${\frak{ g}} \otimes \Bbb C [t, t^{-1}] $.
Faddeev, Reshetikhin and
Takhtajan \cite{FRT2} extended their realization of
$U_q(\frak{g})$ to the quantum loop algebra $U_q(
{\frak{ g}} \otimes \Bbb [t, t^{-1}])$ via a canonical 
 solution  of the Yang-Baxter equation depending on a parameter $z\in\Bbb C$. 
The first realization of the quantum affine algebra $U_q(\hat{
\frak g})$ and its special degeneration called the Yangian were
 obtained by Drinfeld \cite{D2}.
Later Reshetikhin and Semenov-Tian-Shansky
\cite{RS}  incorporated  the central extension in the
previous construction of \cite{FRT1} to obtain  another
realization of the quantum affine algebra $U_q(\hat{ \frak{g}})$.

There is  an automorphism $D_z$ of $U_q(\hat {\frak g})$ defined as 
 $D_z(e_0)= z e_0, D_z(f_0)=z^{-1}f_0, $  and $D_z$ fixes   all other
generators, where $e_o$ and $f_0$ are generators corresponding to the zero node
of its Dynkin diagram.  We also define the map $\Delta_z(a)= (D_z \otimes
id)\Delta(a)$ and $\Delta'_z(a)= (D_z \otimes id)\Delta'(a)$, where 
$a\in U_q(\hat {\frak g})$ and $\Delta'$ denotes the opposite
comultiplication. Let $d$ be an operator such that $d$ commutes with all
 other elements
but has the relations $[d, e_0]=e_0, [d, f_0]=-f_0.  $ 
It is clear that the action of $D_z$ is equivalent to the 
 conjugation by 
$z^d$.  
For the algebra generated by ${U_q(\hat{\frak g})}$ and 
the operator $d$ which we denote by 
$U_q(\tilde {\frak g})$, from the theory of Drinfeld\cite{D2}, 
we know that it has a universal R-matrix $\bar {\frak R}$,
$\bar \frak{ R}$ in $U_q^+(\tilde {\frak g})\hat \otimes  U_q^-(\tilde {\frak g})$
 such that $\bar {\frak R}$ satisfies the properties:
$$\gather
\bar {\frak R}\Delta(a)=\Delta^{\text{op}}(a)\bar {\frak R},\\
(\Delta\otimes \text{id})(\bar {\frak R})=\bar {\frak R}_{13}\bar
 {\frak R}_{23},\tag I.1 \\
(\text{id}\otimes\Delta)(\bar {\frak R})=
\bar {\frak R}_{13}\bar {\frak R}_{12},
  \endgather
$$
where $a\in U_q({\frak g})$,   $\Delta^{op}$ denotes the opposite 
comultiplication, $\bar
{\frak  R}_{12}=\sum_ia_i\otimes b_i\otimes 1=\bar {\frak R}\otimes 1,
\bar {\frak R}_{13}=\sum_ia_i\otimes 1\otimes b_i,  
\bar {\frak R}_{23}=\sum_i1\otimes a_i\otimes b_i=
1\otimes \bar {\frak R}. $
Here $U_q^+(\tilde {\frak g})$ is the subalgebra generated by 
$e_i$, $t_i$ and $d$ and $U_q^-(\tilde {\frak g})$ is the subalgebra
generated by $f_i$, $t_i$ and $d$. 

Let $\bar{\frak R}(z)= (D_z\otimes {\text id})\bar {\frak R}. $ 
Let $C$ be  the central element corresponding to the central
extension of the quantum affine algebra.
 Let  ${\frak R}(z)$=$ q^{-d\otimes C -C\otimes d}\bar
{\frak R}(z)$, then 
we have   ${\frak R}(z)$
$\in$ $U_q(\hat{\frak  b^+})\hat  \otimes 
 U_q(\hat {\frak  b^-})\otimes
{\Bbb C} [[z]]$, such that 
$$ {\frak R}(z) \Delta_z (a)= (D^{-1}_{q^{C_2}} \otimes
D^{-1}_{q^{C_1}}) \Delta'_z (a) {\frak R}(z),$$
$$ (\Delta \otimes I){\frak R}(z)= {\frak R}_{13}(zq^{C_2}){\frak R}_{23}(z),$$
$$(I\otimes\Delta){\frak R}(z)= {\frak R}_{13}(zq^{-C_2}){\frak R}_{12}(z),
\tag I.2 $$
$${\frak R}_{12}(z){\frak R}_{13}(zq^{C_2}/w){\frak R}_{23}(w)=
{\frak R}_{23}(w){\frak R}_{13}(zq^{-C_2}/w){\frak R}_{12}(z).  $$
Here $C_1=C\otimes 1$, $C_2=1\otimes C$,  $U_q(\hat{\frak  b^+})$
is the subalgebra generated by $e_i, t_i$ and  $U_q( {\frak  b^-})$
is the subalgebra generated by $f_i, t_i$. 

Let ${\frak C}= ((D^{-1}_{q^{C_2}} \otimes
D^{-1}_{q^{C_1}}) {\frak R}_{21})
{\frak R}$. Then 
$${\frak C} \Delta (a)=
((D^{-2}_{q^{C_2}} \otimes
D^{-2}_{q^{C_1}}) \Delta (a)){\frak C}.  \tag I.3 $$
We note that $(q^{C_2}\otimes q^{C_1})$
is invariant under the permutation. 
This shows that the action of  
 ${\frak C}=((D^{-1}_{q^{C_2}} \otimes
D^{-1}_{q^{C_1}}) {\frak R}_{21})
{\frak R}$ on a tensor product of two modules 
 is  an  intertwiner which, however,   shifts the  actions of
$e_0$ and $f_0$ by  the constants 
$q^{\mp 2C_2}\otimes q^{\mp 2C_1}$ respectively. 
We should also notice that $D_z\otimes D_z$ acts invariantly on 
${\frak R}$. Let ${\frak R}_{21}(z)=(D_z\otimes 1){\frak R}_{21}.$
Note that ${\frak R}_{21}(z)$ is not equal to $P({\frak R}(z))$, where 
$P$ is the permutation operator. 
 
Let $V$ be a finite dimensional   representation of
$U_q(\hat{\frak g})$.
Let ${\bar L}^{+}(z)= ({\text id} \otimes
\pi_{V})( {\frak R}_{21}(z)),  $
${\bar L}^{-^{-1}}(z)= ({\text id} \otimes
\pi_{V}) {\frak R}_{}(z^{-1}), $
${\frak L}^{+}(z)= (\pi_{V} \otimes{\text id})
( {\frak R}^{-1}_{}(z)), $
${\frak L}^{-^{-1}}(z)= (\pi_V \otimes 
{\text id}) {\frak R}^{-1}_{21}(z^{-1}). $
We know that $
{\bar L}^{-^{-1}}(z) P {\frak L}^{-^{-1}}(z^{-1})=1, $ where
$P$ is the permutation operator. 

$U_q(\hat{\frak g})$ as an algebra 
 is  generated  by operator entries of
${\bar L}^{+}(z)$ and
${\bar L}^{-^{-1}}(z)$,  and it is also  generated by operator entries 
of  ${\frak L}^{+}(z)$ and
${\frak L}^{-^{-1}}(z)$. 
${\bar L}^\pm(z)$ are   used in \cite{FR} to obtained 
q-KZ equation. However,  ${\bar L}^\pm(z)$ and ${\bar\frak { L}}^\pm(z)$
are basically the same. 

 Let 
$${\bar L}({z})= ({\text id} \otimes \pi_{V})((1\otimes
D^{-1}_{q^{C}})
 {\frak R}_{21}({z})){\frak R}(z^{-1}), $$
$$({D^{}_z \otimes 1}) {\bar L}({z})=\bar  L = ({\text id} \otimes
\pi_{V})((1\otimes D^{-1}_{q^{C}} ) {\frak R}_{21}){\frak R},  $$
$${\frak L}({z})= (\pi_{V} \otimes{\text id})
(( D^{}_{q^{C}} \otimes 1 ) {\frak R}^{-1}
({z})){\frak R}^{-1}_{21}(z^{-1}), $$
$$(1\otimes D^{}_z) {\frak L}({z})= {\frak L} =(\pi_{V} 
\otimes{\text id})
((1\otimes D^{}_{q^{C}} ) {\frak R}^{-1}){\frak R^{-1}_{21}}.  $$ 
Then 
$$\eqalign{
\bar R(\frac{z}{w}) {\frak L}^{\pm}_1(z)^{-1}
 {\frak L}^{\pm}_2(w)^{-1} &= {\frak L}^{\pm}_2(w)^{-1} 
{\frak L}^{\pm}_1(z)^{-1}
\bar R(\frac{z}{w}),\cr
\bar R(\frac{zq^{-C}}{w}) {\frak L}^+_1(z)^{-1} {\frak L}^-_2(w)^{-1} &=
 {\frak L}^-_2(w)^{-1} {\frak L}^+_1(z)^{-1} 
\bar R(\frac{zq^{C}}{w})
,} \tag I.4 $$
$${\frak L}_1(z)\bar R(zq^{-2C}/w)^{}{\frak L}_2(w)
\bar R(z/w)^{-1}=  \bar R(z/w)
{\frak L}_2(w)\bar R(zq^{2C}/w)^{-1}{\frak L}_1(z),$$
$${\frak L}({z})(id \otimes \pi_V)\Delta(a)= (D^{}_{q^{2C}}\otimes 1
\Delta(a)){\frak L}(z), $$  
Here $\bar R(z/w)$ is the image  of  ${\frak R}(z/w)$ on $V\otimes
V$. 

We name ${\bar L}(z)$ and  ${\frak L}(z)$ universal Casimir operator
and inverse universal Casimir operator of the quantum algebra  respectively.

The construction of 
a  representation  $U_q({\hat{\frak g}})$  is equivalent to finding 
 a   specific realization  of  the 
operator  ${\bar L}(z)$ or ${\frak L}(z)$, which plays the same role as $L$
 in the case of $U_q({\frak g})$ in \cite{DF2}.
Naturally we would like to find a way to build this ${\bar L}(z)$
or ${\frak L}(z)$ out of the 
intertwiners, just as in the case of spinor and oscillator
representations
of the  quantum groups of types $A$, $B$, $C$ and $D$ in \cite{DF2}. 
For this paper, we will only use ${\frak L}(z)$. 

Let $V_{\lambda,k}$ and $V_{\lambda_1,k}$ 
be two  highest weight representations  of 
$U_q(\hat{\frak g})$ with highest weight $\lambda$ and 
$\lambda_1$ and the center $C$
acting  as a  multiplication by  a number $k$. 
 Let $\Phi$ be an intertwiner:
$\Phi:V_{\lambda_1,k} \longrightarrow V\otimes V_{\lambda,k} ,$
$\Phi (x)=E_1\otimes  \Phi_1(x)+...+ E_n\otimes \Phi_n(x), $
where $x\in V_{\lambda_1, k}$ and $\{ E_i\}$ is the  basis for $V$.
Let $\Phi^*$ be an  intertwiner:
$\Phi^*:V_{\lambda,k} \longrightarrow ^*V\otimes  V_{\lambda_1,k} ,$
$\Phi^* (x)= E_1^*\otimes\Phi_1^*(x)+...+E_n^*\otimes  \Phi_n^*(x)
, $
where $^*V$ is the right dual representation of $V$ of $\hat U_q{({\frak
g}})$, $x\in V_{\lambda,k}$ and $E_i^*$ is the  basis for $^*V$.
 By the right dual representation
of $V$ of $\hat U_q{({\frak
g}})$, we mean the action of $\hat{U_q({\frak
g}})$ on the dual space given by $<av',v>=<v',S^{-1}(a)v>$, for
$a\in \hat U_q{({\frak g}})$, $v\in V$ and $v' \in$ $ ^*V$.
$(1\otimes \Phi)\Phi^*=
\Sigma \Phi_i \Phi^*_j\otimes E^*_j\otimes
E_i$ gives a map
 $\Phi:V_{\lambda,k} \longrightarrow  ^*V \otimes V \otimes 
V_{\lambda,k} .  $

Let us identify $^*V\otimes V$ with End$(V)$ 
by the  map  which maps  the first two components of 
 $V\otimes ^*V\otimes V$ to $\Bbb C$ and fix the last component.
Let $\tilde {\frak L}\in$
 End$(V)\otimes$  End$(V_{\lambda,k})$ be:
$\tilde {\frak L}=
 (\tilde {\frak L}_{ij})= ((D_{q^{2k}}\Phi_j)\Phi_i^*  ), 
  $ where $(D_{q^{2k}}\Phi _i)$ means shifting  the evaluation
representation by constant $q^{2k}$ and we assume $\tilde {\frak L}$ is 
well defined. Then 
 $$\tilde {\frak L} \Delta (a)=  (( D_{q^{2k}}\otimes 1) \Delta(a)) \tilde
{\frak L}. \tag I.5 $$
The key idea  is to identify   $\tilde L$ with  
 $L$ or $\tilde {\frak L}$ with ${\frak L}$
 to obtain representations out of intertwiners, which is how we 
obtained the spinor representations for $U_q(\hat \frak{gl}(n))$ \cite{Di}.

This paper  is arranged as follows. We prepare 
the basic results about two point correlation functions based on 
the solutions of q-KZ equations\cite{DvO} in Section 1. 
Then in  Section 2,  we will 
present the commutation relations of the intertwiners,
identify them as quantum affine Clifford algebras,  and 
finally derive spinor representations. 
We will also discuss its connection with other problems. 

{\bf 1.  R-matrices and correlation functions} 

In this section, we will give the results on the 
two point functions of type II intertwiners, which basically 
come from the results of  
Davies and Okado
\cite{DvO} about the explicit solutions of q-KZ equations of 
the quantum groups of $B_n^{(1)}$ and $D_n^{(1)}$. We will follow the 
notations  in \cite{DvO}. 

Let ${\frak h}^*$, $\Lambda_i,h_i={\alpha}_i^{\vee},\alpha_i,
 \delta=\sum_{i=0}^{n}a_i\alpha_i$
and $d$ be the standard notations in \cite{K}.
The dual Coxeter number is defined as $h^\vee=\sum_{i=0}^n
{a}_i^{\vee}$, where ${a}_i^\vee=a_i$ 
except  for $ B_n^{(1)}$, ${a}_n^{\vee}=a_n/2$,
so that we have $h^{\vee }=2n-1$ for $B_n^{(1)}$ and  $2n-2$ for 
$D_n^{(1)}$. The invariant bilinear form $(\,|\,)$ in \cite{K} is normalized 
so that $(\theta|\theta)=2$,
where $\theta=\delta-\alpha_0$. We put $\rho=\sum_{i=0}^n\Lambda_i$.
For $\lambda\in{ \frak h}^*$ $\bar {\lambda}$ denotes the orthogonal 
projection of $\lambda$ on
$\bar{ {\frak h}}^*$, where  $\bar {{ \frak h}}^*$ is the linear 
span of the classical roots
$\alpha_1,\cdots,\alpha_n$. We introduce an orthonormal
basis $\{\epsilon_1,\cdots,\epsilon_n\}$ of $\bar{{ \frak h}}^*$, 
by which $\alpha_i,
\bar {\Lambda}_i,\bar {\rho}$ are represented below.
$$\eqalign{
\alpha_i&=\epsilon_i-\epsilon_{i+1}\hskip2.1cm(1\le i\le n-1),\cr
&=\epsilon_n\hskip3.06cm(i=n\text{ for }B_n^{(1)}),\cr
&=\epsilon_{n-1}+\epsilon_n\hskip1.95cm(i=n\text{ for }D_n^{(1)}),\cr
\bar{\Lambda}_i&=\epsilon_1+\cdots+\epsilon_i
\hskip1.6cm(1\le i\le n-1\text{ for }B_n^{(1)},
1\le i\le n-2\text{ for }D_n^{(1)}), \cr
&={\epsilon_1+\cdots+\epsilon_{n-1}-\epsilon_n \over2}\quad(i=n-1\text{
 for }D_n^{(1)}),\cr
&={\epsilon_1+\cdots+\epsilon_{n-1}+\epsilon_n \over2}\quad(i=n),\cr
2\bar{\rho}&=(2n-1)\epsilon_1+(2n-3)
\epsilon_2+\cdots+3\epsilon_{n-1}+\epsilon_n
\quad(\text{ for }B_n^{(1)}),\cr
&=(2n-2)\epsilon_1+(2n-4)\epsilon_2+\cdots+2\epsilon_{n-1}\quad(\text{ for
 }D_n^{(1)}).\cr}\tag 1.1 $$

\proclaim{Definition 1.1}
The quantum affine algebra $U_q(\tilde{\frak g})$ is defined as the {\bf
Q}$(q)$--algebra
with 1 generated by $e_i,f_i,t_i(i=0,\cdots,n),q^d$ satisfying
$$
\eqalign{
&[t,t']=0\hskip2.1cm\hbox{for}\quad t,t'\in\{t_0,\cdots,t_n,q^d\},\cr
&t_ie_jt_i^{-1}=q^{(\alpha_i|\alpha_j)}e_j,\quad
t_if_jt_i^{-1}=q^{-(\alpha_i|\alpha_j)}f_j,\cr
&q^de_jq^{-d}=q^{\delta_{j0}}e_j,\,\qquad
q^df_jq^{-d}=q^{-\delta_{j0}}f_j,\cr
&[e_i,f_j]=\delta_{ij}(t_i-t_i^{-1})/(q_i-q_i^{-1}),
\cr
&\sum_{k=0}^b(-)^ke_i^{(k)}e_je_i^{(b-k)}
=\sum_{k=0}^b(-)^kf_i^{(k)}f_jf_i^{(b-k)}=0\quad(i\ne j),\cr}\tag 1.2
$$
where $b=1-<h_i,\alpha_j>$. Here we have set $q_i=q^{(\alpha_i|\alpha_i)/2},
[m]_i=(q_i^m-q_i^{-m})/(q_i-q_i^{-1}),[k]_i!=\prod_{m=1}^k[m]_i,e_i^{(k)}=
e_i^k/[k]_i!,f_i^{(k)}=f_i^k/[k]_i!$.
We denote by $U_q(\hat{\frak g})$ the subalgebra of
$U_q(\tilde{\frak g})$ generated by the elements $e_i,f_i,t_i$ ($i=0,\cdots,n$).
Let $x_i$ be any of $e_i,f_i,t_i$. 
The coproduct $\Delta$ and the antipode $a$ are defined as follows.
$$
\eqalignno{
\Delta(e_i)&=e_i\otimes 1+t_i\otimes e_i,\qquad
\Delta(f_i)=f_i\otimes t_i^{-1}+1\otimes f_i,&  \cr
\Delta(t_i)&=t_i\otimes t_i,\hskip2.05cm \Delta(q^d)=q^d\otimes q^d,\cr
a(e_i)&=-t_i^{-1}e_i,\qquad a(f_i)=-f_it_i,
\qquad a(t_i)=t_i^{-1},\qquad a(q^d)=q^{-d}.&\cr }
$$
\endproclaim

For a representation $(\pi,V)$ of $U_q(\hat{\frak g})$, we put 
$V_z=V\otimes${\bf Q}$(q)[z,z^{-1}]$, and lift $\pi$ to a 
representation $(\pi_z,V_z)$ of $U_q(\tilde{\frak g})$
as follows:
$
\pi_z(e_i)(v\otimes z^n)=\pi(e_i)v\otimes z^{n+\delta_{i0}},
\pi_z(f_i)(v\otimes z^n)=\pi(f_i)v\otimes z^{n-\delta_{i0}},
\pi_z(t_i)(v\otimes z^n)=\pi(t_i)v\otimes z^n, 
\pi_z(q^d)(v\otimes z^n)=v\otimes (qz)^n.
$  Define an index set $J$ by
$J=\{0,\pm1,\cdots,\pm n\},\phantom{B_n^{(1)}}
\hbox{for }{\frak g}=B_n^{(1)};$ 
 $=\{\pm1,\cdots,\pm n\},\phantom{D_n^{(1)}}
\hbox{for }{\frak g}=D_n^{(1)},
$
and set $N=|J|$. We introduce a linear order $\prec$ in $J$ by
$
1\prec2\prec\cdots\prec n\;(\prec0)\prec-n\prec\cdots\prec-2\prec-1.
$
We shall also use the usual order $<$ in $J$.

We shall define the vector representation $(\pi^{(1)},V^{(1)})$ of
$U_q(\hat{\frak g})$.
It is the fundamental representation associated with the first node
in the Dynkin diagram.
The base vectors of $V^{(1)}$ are given by $\{v_j\mid j\in J\}$.
The weight of $v_j$ is given by $\epsilon_j(j>0)$, $-\epsilon_{-j}(j<0)$,
$0\,(j=0)$.  We take $v_1$ as a reference vector.
Set
$s=(-)^n$  for ${\frak g}=B_n^{(1)}$, and 
         $s= (-)^{n-1}$ for ${\frak g}=D_n^{(1)}$.

Denoting the matrix units by $E_{ij}$ {\it i.e.} $E_{ij}v_k=\delta_{jk}v_i$,
the actions of the generators read as follows
($\pi^{(1)}(f_i)=\pi^{(1)}(e_i)^t$):
$$
\eqalignno{
\pi^{(1)}(e_0)&=s(E_{-1,2}-E_{-2,1}),\cr
\pi^{(1)}(t_0)&=\sum_{j\in
J}q^{-\delta_{j1}-\delta_{j2}+\delta_{j,-1}+\delta_{j,-2}}
          E_{jj},\cr
\pi^{(1)}(e_i)&=E_{i,i+1}-E_{-i-1,-i}\quad(1\le i\le n-1),\cr
\pi^{(1)}(t_i)&=\sum_{j\in J}q^{\delta_{ji}-\delta_{j,i+1}+\delta_{j,-i-1}
          -\delta_{j,-i}}E_{jj}\quad(1\le i\le n-1),\cr
\pi^{(1)}(e_n)&=\sqrt{[2]_n}(E_{n0}-E_{0,-n})
          &\hbox{for }{\frak g}=B_n^{(1)},\hskip1.5cm  \cr
        &=E_{n-1,-n}-E_{n,-n+1}
          &\hbox{for }{\frak g}=D_n^{(1)},\hskip1.5cm\cr
\pi^{(1)}(t_n)&=\sum_{j\in J}q^{\delta_{j,n}-\delta_{j,-n}}E_{jj}
          &\hbox{for }{\frak g}=B_n^{(1)},\hskip1.5cm\cr
        &=\sum_{j\in J}q^{\delta_{j,n-1}+\delta_{jn}-\delta_{j,-n}
          -\delta_{j,-n+1}}E_{jj}
          &\hbox{for }{\frak g}=D_n^{(1)}.\hskip1.5cm\cr}
$$

Let $\{v^*_j\mid j\in J\}$ be the canonical dual basis of $\{v_j\mid j\in J\}$.
Then we have the following isomorphism of $U_q(\tilde{\frak g})$--modules.
$$
\eqalignno{
C^{(1)}_\pm:\quad V^{(1)}_{z\xi^{\mp 1}}
&{\buildrel \sim\over\longrightarrow}\bigl(V^{(1)}_z\bigr)^{*a^{\pm1}};\qquad
v_j\mapsto q^{\overline{\pm j}}v_{-j}^*.
&\cr} 
$$
Here $\overline{j}$ is defined by
$\overline{j}=j$, if $j>0$; or $n$, if $j=0$; or $j+N$, if $j<0$.

We also define a matrix $\bar R^{(1,1)}(z)$ on $V^{(1)}\otimes V^{(1)}$ as:
 $$
\eqalign{
\bar R^{(1,1)}(z)=&\sum_{i\ne0}E_{ii}\otimes E_{ii}
      +{q(1-z)\over1-q^2z}\sum_{i\ne\pm j}E_{ii}\otimes E_{jj}\cr
     &+{1-q^2\over1-q^2z}\Bigl(\sum_{i\prec j,i\ne-j}+z\sum_{i\succ j,i\ne-j}
      \Bigr)E_{ij}\otimes E_{ji}  \cr
     &+{1\over(1-q^2z)(1-\xi z)}\sum_{i,j} a_{ij}(z)E_{ij}\otimes E_{-i,-j}.
\cr} \tag 1.3 
$$
Here
$$
\eqalign{
a_{ij}(z)&=(q^2-\xi z)(1-z)+\delta_{i0}(1-q)(q+z)(1-\xi z), (i=j);\cr
         &=(1-q^2)(q^{\overline{j}-\overline{i}}(z-1)+\delta_{i,-j}(1-\xi z))
               ,(i\prec j);\cr
         &=(1-q^2)z(\xi q^{\overline{j}-\overline{i}}(z-1)+\delta_{i,-j}(1-\xi
z)), 
               (i\succ j).\cr} \tag 1.4 
$$

Let $V(\Lambda_i)$ be the irreducible highest weight module with highest weight
$\Lambda_i$. Let $\ket{\Lambda_i}$ be a highest weight vector of
 $V(\Lambda_i)$. We only consider level 1 modules , {\it i.e.} $i=0,1,n$
(and $n-1$ for ${\frak g}=D_n^{(1)}$).
Let $\lambda,\mu$ stand for level 1 weights.
We define intertwiners as the
 the following:
$$
\eqalignno{
\Phit_\lambda^{\mu V^{(1)}}(z)&:\quad
V(\lambda)\longrightarrow V(\mu)\otimes V^{(1)}_z,&\cr
\Phit_\lambda^{V^{(1)}\mu}(z)&:\quad
V(\lambda)\longrightarrow V^{(1)}_z\otimes V(\mu).&\cr}
$$
We call them type I and type II  respectively  depending  on the location 
of $V^{(1)}$.

For the vector representation $V^{(1)}$ we fix the normalisations as
follows:
$$
\eqalign{
\Phit_{\Lambda_1}^{ V^{(1)}\Lambda_0 }(z)\ket{\Lambda_1}
&=v_1\otimes \ket{\Lambda_0}+\cdots\cr
\Phit_{\Lambda_n}^{V^{(1)}\Lambda_n }(z)\ket{\Lambda_n}
&= \alpha^{-1}v_0\otimes \ket{\Lambda_n}+\cdots\cr}
\qquad\eqalign{
&\hbox{for }{\frak g}=B_n^{(1)},D_n^{(1)},\cr
&\hbox{for }{\frak g}=B_n^{(1)},\cr}\tag 1.5
 $$
where
$$
\alpha=\sqrt{[2]_n}.
$$
 $\Phit_{\Lambda_0}^{V^{(1)}\Lambda_1}$
 (${\frak g}=B_n^{(1)} ,D_n^{(1)}$) and
$\Phit_{\Lambda_{n-1}}^{ V^{(1)}\Lambda_n }$,  
$\Phit_{\Lambda_n}^{V^{(1)}\Lambda_{n-1}}$
(${\frak g}=D_n^{(1)}$) are normalised using  Dynkin diagram automorphisms.

As explained in \cite{FF}, it  is necessary to 
shift the grading by 1/2,  in order to 
write down simple commutation relations between those operators.
 We will assume that 
$\tilde \Phi_\mu^{(V^{(1)})_1\lambda}(z)$, for $\mu=\Lambda_{0,1}$, is shifted 
by $z^{\mp1/2}$ respectively, we normalize them such that 
$$
\Phit_{\Lambda_1}^{ V^{(1)}\Lambda_0 }(z)\ket{\Lambda_1}
=z^{-1/2}v_1\otimes \ket{\Lambda_0}+\cdots,   \tag 1.6 $$
$$\Phit_{\Lambda_0}^{ V^{(1)}\Lambda_1 }(z)\ket{\Lambda_0}
=z^{1/2}v_{-1}\otimes \ket{\Lambda_1}+\cdots. $$

With the nomalization, from now on, we will will use $\Phi^{V^{(1)}\mu}_
\lambda$, instead of 
$\tilde \Phi^{V^{(1)}\mu}_\lambda$ to denote those operators. 

Let $\br{\nu|\Phi_\mu^{\nu W_2}(z_2)\Phi_\lambda^{\mu V_1}(z_1)|\lambda}$
be  the expectation value of the following
composition of the intertwiners
$$
V(\lambda)
{\buildrel \Phi_\lambda^{\mu V}(z_1)\over\longrightarrow}
V(\mu)\otimes V_{z_1}
{\buildrel \Phi_\mu^{\nu W}(z_1)\otimes 1\over\longrightarrow}
V(\nu)\otimes W_{z_2}\otimes V_{z_1}
{\buildrel 1\otimes P\over\longrightarrow}
V(\nu)\otimes V_{z_1}\otimes W_{z_2}.
\tag 1.7  $$

Recalling  the isomorphism $C_\pm^{(1)}$, we define
$$
\Phi_\lambda^{\mu(V^{(1)})^{*a^{\pm1}}}(z)
=(C_\pm^{(1)}\otimes {\text id})\cdot
\Phi_\lambda^{\mu V^{(1)}}(z\xi^{\mp1}).
\tag 1.8  $$

\proclaim{Theorem 1.1} \cite{DvO} \cite{O}
 For any possible combination of weights $(\lambda,\mu)$ 
from $\Lambda_0$, $\Lambda_1$, $\Lambda_{n-1}$,$\Lambda_{n}$ 
$$
\Phi_\mu^{(V^{(1)})_2\lambda}(z_1)\Phi_\lambda^{(V^{(1)}\mu)_1}(z_2)
=P\rho(z_1/z_2)\bar R^{(1,1)}(z_1/z_2)
\Phi_\mu^{(V^{(1)})_2\lambda}(z_2)\Phi_\lambda^{(V^{(1)})_1\mu}(z_1),
\tag 1.9  $$
where
$$
\rho(z)=z
{(z^{-1};\xi^2)_\infty(q^{-2}\xi z;\xi^2)_\infty
(\xi z;\xi^2)_\infty(q^{-2}\xi^2z;\xi^2)_\infty \over
( z;\xi^2)_\infty(q^{-2}\xi^2z;\xi^2)_\infty
(   \xi z^{-1};\xi^2)_\infty(q^{-2}\xi z^{-1};\xi^2)_\infty}, 
\tag 1.10  $$ and $P$ is the permutation operator. 
\endproclaim 

The equality of the commutation relations in the theorem is valid  on the 
level of the correlation functions.

{\bf 2. Spinor representations for $B^{(1)}_n$ and $D^{(1)}_n$.}

In \cite{Di}, we presented 
 a construction of spinor representations of 
$U_q(\hat{\frak gl}(n))$. We will construct spinor representations
in a completely parallel way for $U_q(\hat{\frak o}(n))$, which is 
a quantization of corresponding classical constructions in 
\cite{FF}. We will identify the type II intertwiners 
as generators of quantum affine Clifford algebras and construct 
the inverse universal Casimir operators. 

In this section, we will assume that $|q|<1$. 

Let $\Lambda$ be $V(\Lambda_n)$ for the case of $B_n^{(1)}$ and 
$\Lambda$ be $V(\Lambda_{n-1})\oplus V(\Lambda_n)$  or 
$V(\Lambda_{0})\oplus V(\Lambda_1)$ for the case of 
$D_n^{(1)}$. 

Let $\Phi(z)$ be the intertwiner $\Phi_{\Lambda_n}^{V^{(1)}
\Lambda _n}$ from  $\Lambda$ to  $V^{(1)}_z\otimes \Lambda$ 
for the case of $B_n^{(1)}$ and  and 
$\Phi(z)$ be $\Phi_{\Lambda_{n-1}}^{V^{(1)}
\Lambda _n} \oplus \Phi_{\Lambda_n}^{V^{(1)}\Lambda _{n-1}}(z)$ or
 $\Phi_{\Lambda_{0}}^{V^{(1)}\Lambda _1}(z)
\oplus \Phi_{\Lambda_1}^{V^{(1)} \Lambda _{0}}$
as an intertwiner 
from $\Lambda$ to $ V^{(1)}_z \otimes \Lambda$. 

\proclaim{Theorem  2.1} $\Phi^{}(z)$ satisfies the following commutation 
relations,  
$$ F(z_1/z_2)
\Phi^{2}(z_1)\Phi^{1}(z_2)
=G(z_1/z_2)P \bar R^{(1,1)}(z_1/z_2)
\Phi^{2}(z_2)\Phi^{1}(z_1)+F\delta(z_1\xi/z_2), \tag A$$
when $\Lambda$ is $\Lambda_n$ or $\Lambda_{n-1}\oplus \Lambda_n$; 
$$ F(z_1/z_2)
\Phi^{1}(z_1)\Phi^{2}(z_2)
=G(z_1/z_2)P\bar R^{(1,1)}(z_1/z_2)
\Phi^{2}(z_2)\Phi^{1}(z_1)+(z_1/z_2)^{1/2}F\delta(z_1\xi/z_2), \tag B  $$
when $\Lambda$ is $\Lambda_{0}\oplus \Lambda_1$.
$$
F(z)={ {(\xi^{1} z^{-1};\xi^2)_\infty 
(q^{-2} z^{-1};\xi^2)_\infty}
\over {(z^{-1}\xi^2;\xi^2)_\infty(q^{-2}\xi z^{-1};\xi^2)_\infty}}, \tag 2.1 $$
$$
G(z)= \rho(z)F(z),
$$
$P$ is the permutation operator, 
 the matrix $\bar R^{(1,1)}(z_1/z_2)$ is defined for $B_n^{(1)}$ and 
 $D_n^{(1)}$ respectivel as in Section 1, and $F$ is  a nonzero 
vector in the one dimensional invariant 
subspace under the action of the subalgebra $U_q({\frak g})$ 
generated by  $e_i, f_i$ and $t_i$, $i\neq 0$ respectively. 

\endproclaim

In this theorem, $z_i$ are formal variables and the functions on two 
sides have different expansion directions. Because the infinite formal 
power expansion is involved, we define two 
operators to be equal if the action of these two operators 
 coincide on any vector in $\Lambda$.  

{\bf Proof.}
 The proof is basically the same as in the case of \cite{Di}. 
First,  we look at the matrix coefficients corresponding to the highest 
weight vectors of the two sides of the 
equalities above. With the precise expression  of the correlation functions, 
we can see clearly  that the left hand side is in the form of polymonials of 
$z_1/z_2$ over $(1-z_1\xi/z_2)$. With Theorem 1.1 in the last section, 
it is clear that the first term of the  right hand side is the same 
as  that of the other side, but in a different expansion direction. 
If we move the first term of the right side to the left hand side, we 
can explicitly calculate to show  that this equality is valid.  
So we can prove that this statement is true 
for the coefficients corresponding to the highest weight vectors. 
However, we also know that $F\delta(z_1\xi/z_2)$ and 
$(z_1/z_2)^{1/2}F\delta(z_1\xi/z_2)$  are invariant under the  action 
of the quantum affine algebra,
 due to the fact that $F$ becomes an invariant subspace
when $z_1\xi =z_2$ and $\delta(z_1/z_2)f(z_1,z_2)$=
$\delta(z_1/z_2)f(z_1,z_1)$=$\delta(z_1/z_2)f(z_2,z_2)$ for any 
polymonial $f(z_1, z_2)$. 
Then, we know that $F\delta(z_1\xi/z_2)$ and 
 $(z_1/z_2)^{1/2}F\delta(z_1\xi/z_2)$ are also intertwiners. 
So we can use the quantum affine algebra action to prove that 
this equality is valid for any matrix coefficient. Thus we prove the equality. 

\proclaim{Definition 2.1} 
Quantum affine Clifford algebra of $\Bbb Z$-type
 is an associative algebra generated 
by $\Psi_i(m)$, $i\in J$, $m\in \Bbb Z$, where $J$ is defined in the 
Section 1. Let 
$\Psi(z)=\Sigma_{m}\Psi_i(m)z^{-m}\otimes E_i$, where $E_i$ is a base for 
$V^{(1)}$. Then relations are 
$$ F(z_1/z_2)
\Psi^{}(z_2)\Psi^{}(z_1)
=G(z_1/z_2)\bar R^{(1,1)}(z_1/z_2)P
\Psi^{}(z_1)\Psi^{}(z_2)+F\delta(z_1\xi/z_2) . \tag 2.2 $$
Here $F(z)$ and $G(z)$ are the same as defined above. P is the 
permutation operator and $R^{(1,1)}$ are on $\Bbb C^N$ respectively. 
\endproclaim

\proclaim{Definition 2.2} 
Quantum affine Clifford algebra of $\Bbb Z+1/2$-type
 is an associative algebra generated 
by $\Psi_i(m)$, $i\in J$, $m\in \Bbb Z+1/2$, where $J$ is defined in the 
Section 1 with $N=2n$. Let 
$\Psi(z)=\Sigma_{m}\Psi_i(m)z^{-m}\otimes v_i$. Then relations are 
$$ F(z_1/z_2)
\Psi^{}(z_2)\Psi^{}(z_1)
=G(z_1/z_2)\bar R^{(1,1)} P (z_1/z_2)
\Psi^{}(z_1)\Psi^{}(z_2)+(z_1/z_2)^{1/2}F\delta(z_1\xi/z_2) . \tag 2.3 $$
Here $F(z)$ and $G(z)$ are the same as defined above. P is the 
permutation operator and $R^{(1,1)}$ are on $C^{2n}$ only. 
\endproclaim 

In some sense, 
the algebras defined above 
are not rigorously defined, due to the fact that infinite expansions
are involved. However, those definitions  makes sense, 
if we look at a representation of this algebra on 
a space  generated by one vector such that the operators of positive degrees
are  annihilators of this vector. We can also use the specific base as 
in \cite {DvO} to specify those relations to define an ordered base,
with which we can define an  inverse limit type topology. 

\proclaim{Theorem 2.2}
Quantum affine Clifford algebra of $\Bbb Z$-type
is isomorphic to the algebra generated by $\Phi(z)$ for corresponding
cases of type (A). 
Quantum affine Clifford algebra of $\Bbb Z+1/2$-type
is isomorphic to the algebra generated by $\Phi(z)$ for corresponding
cases of type (B).
\endproclaim 

Clearly what we really have to show here is 
that $\Phi(z)$ gives faithful representations 
of those algebras, or we can say that all the relations between $\Phi_j(n)$
are included in the commutation relations (A) or (B). 

{\bf 
Proof.} The proof is basically the same as in \cite{Di}. We can use the base 
as in \cite{DvO} to derive a specific basis for those quantum Clifford 
algebras, where the key observation is that 
${{z(1-1/z)}\over(1-z)}=-1$. This  ensures a wedge type relations
between the generators that will show that this algebra has 
the wedge type character. By comparing characters, we can show that 
it is a faithful  representation. 

Let $\Phi^*(z)= (C_+^{(1)}\otimes {\text id})\cdot
\Phi(z\xi^{- 1}).  $
The locations  of 
the poles of the 
correlation functions of $\Phi^2(z_1)\Phi^{*1}(z_2) $ clearly 
do not include the line  $z_1q^2=z_2$. From the commutation 
relations and the condition that  $|q|<1$, 
the multiplication of $ \Phi^*(z)$ and
$\Phi(zq^2)$  are  well  defined. 
Let $E_i$, $i\in J$, be the base of $V^{(1)}$  as 
defined in Section 1 and 
$E^*_i$ be its right dual base. 
Thus we can say that 
 $(D^{}_{q^{2}}\Phi_i)\bar \Phi^*_j E^*_j\otimes E_i = 
(1\otimes 1\otimes D_z^{-1})
(1\otimes \Phi(zq^2))\bar \Phi^*(z)$ is well defined. 
Here $\Phi(z)=\Sigma E_i\otimes \Phi_i(z)$ and 
$\Phi^*(z)=\Sigma E^*_i\otimes \Phi_i(z)$.

Let $\tilde {\frak L}(z)$=$(1\otimes \Phi(zq^2))\bar \Phi^*(z)$. 
From \cite{KKMMNN}, we 
know that $V^{(1)} \otimes \Lambda_i$ are irreducible. This shows that 
the dimension of space of operators 
$X: V^{(1)}\otimes  \Lambda \longrightarrow  V_{q^{2}}\otimes 
  \Lambda,$
which  satisfy the relation :
$$ X \Delta(a)= ( D^{2}_q\otimes 1 ) \Delta (a) X, \tag 2.4  $$ 
is $1$ for the case of $B^{(1)}_n$  and $2$ for the case of
$D^{(1)}_n$. Due to our normalization, it turns out that 
for the case of $D^{(1)}_n$, those two constants are the same, which
can   be determined  by looking at their
actions  on the highest weight vectors.

\proclaim{Theorem 2.3}
 $$  {\frak L}(z)=c (D_{z}\otimes 1) \tilde {\frak  L},  \tag 2.5  $$
where $c=\frac {{\text tr}(v_0,(D^{}_{q^{2}} \Phi^*) \Phi v_0)} {
  {\text tr}(v_0, {\frak L}(1)v_0)} $ and  $v_0$ is a highest weight 
vector of $\Lambda$. 
\endproclaim

We can remove the $c$ in the formula above 
 by normalizing  $\Phi(z)$ and $\Phi^*(z)$ again. 

Combining all the results above, it is clear that we can start from 
abstract algebra as in Definition 2.1 and 2.2, define the Fock space,
 which is generated by operators of negative degree
and some zero degree operators, then the action of $\frak {L}$ can be 
derived. This gives the quantized spinor representations, which degenerate
into classical constructions in \cite{FF}. 

In \cite{Di}, we 
 derive  commutation for the quantized 
affine  Clifford algebras coming from the spinor representations
of $U_q(\hat \frak{g})$. In the
 classical case\cite{FF}, they completely coincide. It is interesting to 
compare those two algebras, which we believe are basically the same. 
In \cite{B}, one type of our spinor representations was obtained using vertex 
operator construction. It will be interesting 
 to compare our constructions with that in \cite{B}, which should 
lead to another version of quantum Boson-Fermion correspondence, for 
which the result in \cite{JKK} will be very useful. 
As we explain in the introduction, with the knowledge of corresponding 
q-KZ equations, we expect that all the constructions in \cite{FF} 
can be derived in the same way. On the other hand, the further 
development of our theory should lead to the theory of quantization 
of vertex operator algebras\cite{FLM}, for which our new algebras 
should provide the proper example, which is closely related to the 
theory of formal factors in  massive quantum field theory\cite{Sm}. 
By now, our methods are used to construct representations, which are 
deformation of certain classical constructions. We do not know if it is 
possible to use this type of algebraic  constructions  to obtain 
representations, which are not known  for the classical cases.

{\bf Acknowledgments} 
I would like to  thank Professor M.  Okado for his discussion and
advice, especially for Theorem 1.1, which I should attribute to him.

\end